\documentclass{amsart}
\title{Calibrating the Negative Interpretation}
\author{Joan Rand Moschovakis \tt(draft September 9, 2021)}

\usepackage{amsmath}
\address{Occidental College (Emerita)}
\email{joan.rand@gmail.com}

\begin{document}
\maketitle

\section{What this essay is about}
G\"odel and Gentzen proved by negative translations that classical Peano arithmetic {\bf PA} is equiconsistent with its intuitionistic subsystem, Heyting arithmetic {\bf HA}.  By hereditarily replacing $\mathrm{A \vee B}$ by its classical equivalent $\mathrm{\neg (\neg A ~\&~ \neg B)}$, and $\mathrm{\exists x A(x)}$ by its classical equivalent $\mathrm{\neg \forall x \neg A(x)}$, they showed that the negative fragment of {\bf HA} (with only the logical symbols $\mathrm{\&, \neg, \rightarrow, \forall}$ and their axioms and rules) faithfully interprets {\bf PA} in the following sense: the negative translations of the mathematical axioms of {\bf PA} are provable in {\bf HA}, the classical logical axioms and rules for $\mathrm{\&, \neg, \rightarrow}$ and $\mathrm{\forall}$ are correct by intuitionistic logic, and every formula of the full language is provably equivalent in {\bf PA} to its negative translation.\footnote{G\"odel \cite{go1933} also translated $\mathrm{A \rightarrow B}$ hereditarily by $\mathrm{\neg (A ~\&~ \neg B)}$, but Gentzen \cite{ge1933} did not.  This paper is based on the simpler Gentzen translation, and on Kleene's axiomatization of intuitionistic and classical logic, arithmetic and two-sorted number theory in \cite{kl1952} and \cite{klve1965}.}

G\"odel \cite{go1933} interpreted this result as showing that intuitionistic arithmetic {\em contains} classical arithmetic via his ``somewhat deviant'' interpretation.  He observed that the absence of a corresponding result for intuitionistic and classical theories of numbers and number-theoretic functions results from mathematical and philosophical, rather than logical, differences.
For example, the negative translation
\[\mathrm{\forall x \neg (\neg \forall y \alpha(\langle x,y \rangle) = 0 ~\&~ \neg \neg \forall y \alpha(\langle x,y \rangle) = 0)}\]
of the instance $\mathrm{\forall x(\forall y \alpha(\langle x,y \rangle) = 0  \vee \neg \forall y \alpha(\langle x,y \rangle) = 0)}$ of the law of excluded middle is provable by intuitionistic logic, and $\mathrm{\forall x (A(x) \vee \neg A(x)) \rightarrow \exists \beta \forall x (\beta(x) = 0 \leftrightarrow A(x))}$ is provable in the classically correct part of Brouwer's intuitionistic analysis {\bf I} as formalized and developed by Kleene and Vesley in \cite{klve1965}, but the negative translation
\[\mathrm{\forall \alpha \neg \forall \beta \neg \forall x (\beta(x) = 0 \leftrightarrow \forall y \alpha(\langle x,y \rangle) = 0)}\]
of  $\mathrm{\forall \alpha \exists \beta \forall x (\beta(x) = 0 \leftrightarrow \forall y \alpha(\langle x,y \rangle) = 0)}$
is neither provable nor refutable in {\bf I}.\footnote{For the relative independence of
$\mathrm{\forall \alpha \neg \forall \beta \neg \forall x (\beta(x) = 0 \leftrightarrow \forall y \alpha(\langle x,y \rangle) = 0)}$ from {\bf I} cf. \cite{jrm1971}, \cite{jrm2010}.}

Suppose {\bf S} is a subsystem of Kleene's {\bf I} which (unlike {\bf I} itself) is consistent with classical logic.  Then the question is: exactly what must be added to {\bf S} in order to prove the Gentzen negative interpretations of its axioms, hence of its theorems? The goal is to find a simple characterization of the precise constructive cost of expanding {\bf S} to include a faithful copy of its classical twin {\bf S$^\circ$} $\equiv$ {\bf S} + $\mathrm{(\neg \neg A \rightarrow A)}$.

\subsection{Definitions}A formal system {\bf S} based on intuitionistic logic is {\em classically consistent} if and only if {\bf S} + $\mathrm{(\neg \neg A \rightarrow A)}$ is consistent.  The {\em classical content} $\mathrm{E}^g$ of a formula $\mathrm{E}$ is its Gentzen negative interpretation, and the {\em classical content} $\mathrm{\Gamma}^g$ of a set $\mathrm{\Gamma}$ of formulas is the closure under intuitionistic logic of the set $\mathrm{\{E}^g\mathrm{: E \in \Gamma\}}$.  The {\em minimum classical extension} {\bf S}$^{+g}$ of a classically consistent formal system {\bf S} is the closure under intuitionistic logic of $\mathrm{{\bf S}\cup{\bf S}}^g$.

If {\bf S} is an axiomatic system based on intuitionistic logic and $\mathrm{A_1, \ldots, A_n}$ is a list of formulas and (logical or mathematical) schemata, then {\bf S} + $\mathrm{A_1 + \ldots + A_n}$ is the formal system obtained by adding $\mathrm{A_1, \ldots, A_n}$ to the axioms of {\bf S}.  For easier comprehension, the negative translations  $\mathrm{\neg \forall x \neg}$, $\mathrm{\neg \forall \alpha \neg}$ of existential quantifiers will sometimes be replaced by their intuitionistic equivalents $\mathrm{\neg \neg \exists x}$, $\mathrm{\neg \neg \exists \alpha}$ respectively.

\subsection{The example of intuitionistic analysis}   By viewing the choice sequence variables $\mathrm{\alpha, \beta, \ldots}$ of the language $\mathcal{L}$({\bf I}) of {\bf I} alternatively as variables over classical one-place number-theoretic functions, restricting the language and logic by omitting $\mathrm{\vee}$ and $\mathrm{\exists}$ with their axioms and rules, and replacing each mathematical axiom of a classically consistent subsystem {\bf S} of {\bf I} by its negative translation, one obtains a classically equivalent copy {\bf S$^g$} of {\bf S$^\circ$} within {\bf S$^{+g}$}. In particular, if {\bf B} is the subsystem of {\bf I} which omits the continuous choice axiom schema CC$_{11}$ (``Brouwer's Principle for a Function,'' axiom schema $^x$27.1 of \cite{klve1965}) then {\bf B$^\circ$} $\equiv$ {\bf B} + $\mathrm{(\neg \neg A\rightarrow A)}$ is classical analysis with countable choice, and {\bf B$^{+g}$} contains a negative version {\bf B}$^g$ of {\bf B$^\circ$}.

The goal here differs from Kleene's in \cite{kl1965} where he showed that {\bf I} is consistent with all purely arithmetical formulas, and all negations of prenex formulas, of the full language $\mathcal{L}$({\bf I}) which are provable in {\bf B$^\circ$}. Any system {\bf S} in $\mathcal{L}$({\bf I}) which is based on intuitionistic logic and has a classical $\omega$-model (a model with standard integers) may be called {\em classically sound}. The minimum classical extension {\bf S$^{+g}$} of any subsystem {\bf S} of {\bf I} for which classical Baire space is an $\omega$-model is classically sound, contains only the essential intuitionistically dubious principles, and is consistent with {\bf I}.

Some easy consequences of continuous choice, such as the axiom schema DC$_1$ of dependent choice for sequences and  Troelstra's neighborhood function principle NFP, are true in classical Baire space; so {\bf B} + DC$_1$ and {\bf B} + NFP are classically sound subsystems of {\bf I} apparently extending {\bf B}.  Their minimum classical extensions, which will be partially analyzed in a later section, are also consistent with {\bf I}.

\subsection{Mathematically significant extensions of intuitionistic analysis} An extension of {\bf B} may be consistent with {\bf I} without being a  subsystem of {\bf I}.  If {\em either} of the two axioms about to be described is added to {\bf B}, the classical content does not change because the new axioms are true by classical logic alone.  In each case the result is consistent with {\bf I} by an appropriate realizability interpretation, but if {\em both} are added to {\bf B}, the result (with the same classical content) is {\em inconsistent} with {\bf I}.

Kleene established the consistency of {\bf I} ($\equiv$ {\bf B} + CC$_{11}$) relative to {\bf B} by means of function-realizability (Theorem 9.3(a) of \cite{klve1965}).  The strong form
\[\mathrm{MP_1. \; \; \; \; \forall \alpha (\neg \forall x \neg \alpha(x) = 0 \rightarrow \exists x \alpha(x) = 0)}\]
of Markov's Principle is self-realizing over {\bf B}, hence consistent with {\bf I}.  All negative statements true in classical Baire space are realizable by Lemma 8.4(a) of \cite{klve1965}, thus ({\bf B} + MP$_1$)$^{+g}$ + CC$_{11}$ = {\bf B}$^{+g}$ + MP$_1$ + CC$_{11}$ is consistent relative to {\bf B}$^{+g}$ + MP$_1$.

Brouwer refuted Markov's Principle using a ``creating subject'' argument.  Kleene proved  {\bf I} $\not\vdash$ MP$_1$ in \cite{klve1965}, using a typed modification of function-realizability he called ``$_S$realizability.''  Vesley \cite{ve1970} proposed adding to {\bf I} the schema
\begin{multline*}
\mathrm{VS. \; \; \forall \alpha \forall x \exists \beta (\overline{\beta}(x) = \overline{\alpha}(x) ~\&~ \neg A(\beta))}\\
\mathrm{\rightarrow [\forall \alpha (\neg A(\alpha) \rightarrow \exists \beta B(\alpha, \beta)) \rightarrow \forall \alpha \exists \beta (\neg A(\alpha) \rightarrow B(\alpha, \beta))]}
\end{multline*}
where $\mathrm{\overline{\alpha}(x)}$ codes the first $\mathrm{x}$ values of the function $\mathrm{\alpha}$, and $\mathrm{\beta}$ is not free in $\mathrm{\neg A(\alpha)}$.  VS is $_S$realizable, hence consistent with {\bf I}, and {\bf I} + VS $\vdash$ $\mathrm{\neg MP_1}$. Negative statements true in classical Baire space are $_S$realizable by Lemma 10.7 of \cite{klve1965}, so ({\bf B} + VS)$^{+g}$ + CC$_{11}$ = {\bf B$^{+g}$} + VS + CC$_{11}$ is consistent and refutes $\mathrm{MP_1}$.

These examples illustrate the mathematical freedom gained by separating the constructive language from the classical language, entirely eliminating the need for classical logic.  The modification GC$_1^\mathrm{neg}$ of Troelstra's principle GC$_1$ of generalized continuous choice\footnote{GC$_1^\mathrm{neg}$ simply replaces ``almost negative'' by ``negative'' in the statement of GC$_1$, cf.\cite{tr1973}.} characterizes Kleene's function-realizability over {\bf B} + MP$_1$, so {\bf B} + MP$_1$ + GC$_1^\mathrm{neg}$ is a consistent extension of {\bf I} + MP$_1$.

{\bf B} + GC$_1^\mathrm{neg}$ proves that every partial functional which is defined at least on a negative dense subspecies of the intuitionistic continuum (e.g. on all sequences which are not eventually monotone) has a continuous partial extension.  In contrast, {\bf I} + VS $\equiv$ {\bf B} + VS + CC$_{11}$ proves the stronger result that every partial functional defined at least on a negative dense subspecies of the intuitionistic continuum has a continuous total extension, although GC$_1^\mathrm{neg}$ is stronger than CC$_{11}$ and {\bf B} + VS has the same classical content as {\bf B}.

\subsection{Additional examples and related work}  A basic axiomatization of the recursive sequences {\bf MRA}, and its minimum classical extension, are studied in this article.  Intuitionistic arithmetic of arbitrary finite types {\bf HA$^\omega$}, Troelstra's {\bf EL}, Bishop's constructive analysis, and three versions of Brouwer's bar theorem in the context of {\bf B} and {\bf I} are discussed in \cite{jrmgvf2021}.  Vafeiadou's results in that article show that minimum classical extensions of consistent but classically unsound theories like {\bf I} should be maximally consistent for the negative language.  Classically sound extensions of {\bf B} which are subsystems of {\bf I} or consistent with {\bf I}, and classically sound theories such as {\bf MRA} which are inconsistent with {\bf B}, have more interesting minimum classical extensions.

A seminal analysis of double negation shift and the negative interpretation of countable choice, in the context of {\bf HA$^\omega$}, was carried out by  Berardi, Bezem and Coquand in \cite{bebeco1998}.  The recent, technical \cite{fuko2018} treats weak nonconstructive principles in the context of {\bf EL}, {\bf HA} or {\bf HA$^\omega$}.  The bibliographies of both point to related work.  For a precise comparison of Troelstra's {\bf EL} and other weak versions of intuitionistic analysis with the systems treated here see \cite{gvfms}, \cite{gvf2012}.  In \cite{lo2009} I. Loeb analyzes a consequence of VS from a constructive reverse mathematics perspective.

\section{The constructive core of intuitionistic analysis}
Like Brouwer, Bishop worked informally, but it seems unlikely that he would have objected to the mathematical content of any of the axioms or axiom schemas of Kleene's neutral basic system {\bf B} except the principle of bar induction. Bishop used countable choice routinely, so Kleene's strongest countable choice axiom schema ($^x$2.1 in \cite{klve1965}):
\[\mathrm{AC_{01}. \; \; \;  \forall x \exists \alpha A(x,\alpha) \rightarrow \exists \beta \forall x A(x,\lambda y.\beta(\langle x,y \rangle))}\]
may be assumed to hold in constructive analysis, with its consequence ($^*$2.2 in \cite{klve1965}):
\[\mathrm{AC_{00}. \; \; \;  \forall x \exists y A(x,y) \rightarrow \exists \alpha \forall x A(x,\alpha(x))}\]
for all formulas $\mathrm{A(x,\alpha)}$ and $\mathrm{A(x,y)}$ of the language, with free variables of both types allowed and with the appropriate conditions on the distinguished variables (e.g. for AC$_{00}$: $\mathrm{\alpha, x}$ must be free for $\mathrm{y}$ in $\mathrm{A(x,y)}$).

Weaker subsystems of {\bf B} are distinguished by restrictions on AC$_{00}$, which in turn determine the classical omega-models of the subsystems.  Classical omega-models are important for constructive analysis because (a) Bishop's work is consistent with classical  mathematics, and (b) the simplest assumption is that the constructive natural numbers are standard.

\subsection{Two-sorted intuitionistic arithmetic {\bf IA$_1$}}
This is an extension of Kleene's first-order intuitionistic arithmetic {\bf IA$_0$} (\cite{kl1952} p. 82).  {\bf IA$_1$} adds variables $\mathrm{\alpha, \beta, \gamma, \ldots}$ over one-place number-theoretic functions, quantifiers $\mathrm{\forall \alpha, \exists \alpha}$ with their (intuitionistic) logical axioms and rules, and finitely many constants for primitive recursive function(al)s with their defining axioms. Additional primitive recursive function constants, with their definitions, may be added as needed.

Terms (of type 0) and functors (of type 1) are defined inductively.  Church's lambda symbol makes it possible to define primitive recursive functors from terms.  There is an axiom schema of lambda-reduction $\mathrm{(\lambda x.t(x))(s) = t(s)}$ (where $\mathrm{t(x), s}$ are terms, and $\mathrm{s}$ is free for $\mathrm{x}$ in $\mathrm{t(x)}$).

Equality at type 0 is a primitive notion, and is decidable in {\bf IA$_1$}.  Equality at type 1 is defined extensionally by  $\mathrm{\alpha = \beta \equiv \forall x (\alpha(x) = \beta(x))}$, and {\bf IA$_1$} includes the open equality axiom $\mathrm{\forall x \forall y (x = y \rightarrow \alpha(x) = \alpha(y))}$.\footnote{{\bf IA$_1$} is the ``least subsystem'' {\bf L} of {\bf I} in \cite{jrmphd}, \cite{kl1969}.  It is defined precisely in \cite{gvf2012}, \cite{jrmgvf2012}.}

The primitive recursive infinite sequences form a classical omega-model of {\bf IA$_1$}.

\subsection{Intuitionistic recursive analysis IRA} Vafeiadou proved in (\cite{gvf2012}) that Troelstra's formal system {\bf EL} (\cite{tr1973}, \cite{tvd1988}) of elementary constructive analysis and the subsystem  {\bf IRA} $\equiv$ {\bf IA$_1$} + QF-AC$_{00}$ of Kleene's {\bf B} have a common definitional extension, where QF-AC$_{00}$ (``quantifier-free countable choice'') restricts AC$_{00}$ to formulas $\mathrm{A(x,y)}$ containing no sequence quantifiers, and only bounded number quantifiers. {\bf IRA} can also be axiomatized by adding to {\bf IA$_1$} a single axiom, either
\[\mathrm{\forall \rho [\forall x \exists y \, \rho(\langle x,y \rangle) = 0 \rightarrow \exists \alpha \forall x \, \rho(\langle x,\alpha(x) \rangle) = 0]} \; \; \;\ \mbox{or}\]
\[\mathrm{\forall \rho [\forall x \exists y \, \rho(\langle x,y \rangle) = 0 \rightarrow \exists \alpha \forall x [\rho(\langle x,\alpha(x) \rangle) = 0 ~\&~ \forall z < \alpha(x) \, \rho(\langle x,z \rangle) \neq 0]],}\]
asserting that the universe of sequences is closed under unbounded constructive search.\footnote{Veldman prefers the unbounded search axiom to the schema QF-AC$_{00}$ for his system {\bf BIM} of intuitionistic recursive analysis (cf. \cite{vel2014}).}

The general recursive infinite sequences provide a natural classical omega-model of intuitionistic recursive analysis {\bf IRA}.

\subsection{Countable comprehension and arithmetical countable choice}
Stronger than QF-AC$_{00}$ over {\bf IA$_1$}, but weaker than AC$_{00}$, is {\em countable comprehension} or ``unique choice''
\[\mathrm{AC_{00}!. \;  \; \; \forall x \exists ! y A(x,y) \rightarrow \exists \alpha \forall x A(x,\alpha(x)),}\]
where $\mathrm{\exists ! y A(x,y)}$ always abbreviates $\mathrm{\exists y A(x,y) ~\&~ \forall y \forall z (A(x,y) ~\&~ A(x,z) \rightarrow y = z)}$.
Since quantifier-free formulas are decidable in {\bf IA$_1$}, the hypothesis of an instance of QF-AC$_{00}$ provides unique least witnesses for the corresponding instance of AC$_{00}$! and so AC$_{00}$! entails QF-AC$_{00}$ -- but not conversely.

Vafeiadou (\cite{gvf2012}, \cite{gvfms}) proved that AC$_{00}!$ is equivalent over {\bf IRA} to the schema
\[\mathrm{CF_d. \; \; \; \forall x (A(x) \vee \neg A(x)) \rightarrow \exists \alpha \forall x [\alpha(x) \leq 1 ~\&~ (\alpha(x) = 0 \leftrightarrow A(x))],}\]
asserting that every analytically definable subset of the natural numbers with a decidable membership relation has a characteristic function.  The converse of $\mathrm{CF_d}$ is provable in {\bf IA$_1$}.

It follows that {\bf IA$_1$} + AC$_{00}$! and {\bf IA$_1$} + AC$_{00}$ have the same classical omega-models, including all analytically definable infinite sequences.

A formula of the two-sorted language is called {\em arithmetical} if it contains only number quantifiers; free variables of both types are permitted. The {\em arithmetical countable choice} schema AC$_{00}^{Ar}$ restricts AC$_{00}$ to arithmetical formulas $\mathrm{A(x,y)}$, and {\em arithmetical comprehension} AC$_{00}^{Ar}$! is the corresponding restriction of AC$_{00}$!.

The arithmetical sequences provide a classical omega-model of {\bf IA$_1$} + AC$_{00}^{Ar}$ (and of {\bf IA$_1$} + AC$_{00}^{Ar}$!).

\subsection{Full countable choice and function comprehension} The schema AC$_{01}$ expresses countable choice for functions. AC$_{01}$! (with $\mathrm{\forall x \exists ! \alpha A(\overline{\alpha}(x))}$ as hypothesis) expresses the corresponding function comprehension principle, where in general $\mathrm{\exists ! \alpha B(x)}$ $\equiv$ $\mathrm{\exists \alpha B(x) ~\&~ \forall \alpha \forall \beta (B(\alpha) ~\&~ B(\beta) \rightarrow \forall x \alpha(x) = \beta(x))}$.

While AC$_{00}$ is weaker than AC$_{01}$ both classically and intuitionistically, AC$_{00}$! is equivalent to AC$_{01}$! over {\bf IA$_1$}.\footnote{cf. \cite{jrmphd}, \cite{jrm1967} where {\bf M} = {\bf IA$_1$} + AC$_{00}$! is proposed as a minimal base theory for constructive analysis.  However, Troelstra \cite{tr1973} observed that Kleene's formalization \cite{kl1969} of the theory of recursive functionals in {\bf M} could equally well be done in {\bf EL}, hence in {\bf IRA}.} Although Kleene chose AC$_{01}$ as an axiom schema for {\bf B}, he observed in \cite{klve1965} that in all but one instance AC$_{00}$ would have sufficed.  It could be interesting to look for essential uses of the stronger principle in constructive and intuitionistic mathematics.

AC$_{00}$ is equivalent over {\bf IRA} to dependent choice for numbers
\[\mathrm{DC_0.  \; \; \; \forall x \exists y A(x,y) \rightarrow \forall x \exists \alpha (\alpha(0) = x ~\&~ \forall y A(\alpha(y), \alpha(y+1))).}\] Over {\bf IA$_1$} + AC$_{00}$ + $\mathrm{(\neg \neg A \rightarrow A)}$, DC$_0$ is equivalent to classical bar induction BI$^\circ$ (to be described in the next section) using $^\ast$26.1$^\circ$ in \cite{klve1965} together with a straightforward converse argument.

It follows that every classical $\omega$-model of {\bf IA$_1$} + AC$_{01}$ is also an $\omega$-model of {\bf B}, since {\bf IA$_1$} + AC$_{01}$ $\vdash$ AC$_{00}$. Moreover, {\bf IA$_1$} + AC$_{00}$, {\bf IA$_1$} + AC$_{00}!$ and {\bf IA$_1$} + BI$_\mathrm{d}$ all have the same classical $\omega$-models, where BI$_\mathrm{d}$ is intuitionistic bar induction with a decidable bar (to be described in the next section).

\section{Brouwer's principles of bar and fan induction}
In addition to full mathematical induction and the principle of countable choice, Brouwer believed he could justify another classically sound principle known as the ``bar theorem.'' Kleene analyzed Brouwer's proof of this principle and found it to be circular.  Kleene's {\bf B} has an axiom schema of bar induction in four versions, which are equivalent over {\bf IA$_1$} + AC$_{00}$!.  Each has the general form\footnote{In Kleene's primitive recursive coding $\mathrm{\langle a_0,\ldots,a_n \rangle = \Pi_{j=0}^{j=n} p_j^{a_j}}$ where $\mathrm{p_j}$ is the j$th$ prime, and $\mathrm{(\langle a_0,\ldots,a_n \rangle)_j = a_j}$.  ``Sequence numbers'' w satisfying $\mathrm{Seq(w) \equiv \forall j < lh(w) \, (w)_j \neq 0}$ uniquely code finite sequences of numbers, where
$\mathrm{lh(w) = \Sigma_{j < w} sg((w)_j)}$ and $\mathrm{sg(n) = 1 \stackrel{.}{-}(1 \stackrel{.}{-} n)}$.  1 codes the empty sequence, $\mathrm{\langle a_0+1,\ldots,a_n+1 \rangle}$ codes $\mathrm{(a_0,\ldots,a_n)}$ and $\mathrm{*}$ denotes concatenation. $\mathrm{\overline{\alpha}(0) = 1}$ and $\mathrm{\overline{\alpha}(n+1) = \langle \alpha(0)+1, \ldots, \alpha(n)+1 \rangle}$.}
\begin{multline*}
\mathrm{BI. \; \; \; \forall \alpha \exists x R(\overline{\alpha}(x)) ~\&~ \forall w (Seq(w) ~\&~ R(w) \rightarrow A(w))} \\
\mathrm{ ~\&~ \forall w (Seq(w) ~\&~ \forall s A(w \ast \langle s+1 \rangle) \rightarrow A(w)) \rightarrow A(1),}
\end{multline*}
where $\mathrm{R(w)}$ is the basis (or bar) predicate and $\mathrm{A(w)}$ is the inductive predicate.\footnote{Later Kreisel and Troelstra \cite{krtr1970} developed a competing formal system for Brouwer's analysis in which the ``bar theorem'' was treated as a principle of generalized inductive definition; cf. \cite{ka2019a}.} As usual, free variables of both types are allowed.

{\em Classical bar induction} BI$^\circ$ places no restrictions on $\mathrm{R(w)}$. Kleene observed that BI$^\circ$ conflicts with Brouwer's continuity principle so some restriction is necessary in the intuitionistic context.

Brouwer used bar induction to prove his ``fan theorem,'' which (together with the assumption that every full function is pointwise continuous) allowed him to conclude that every function completely defined on the closed unit interval is uniformly continuous there.  The {\em full fan theorem} (\cite{klve1965} $^\ast$27.9), which is provable in {\bf I} for all predicates $\mathrm{R(w)}$ in which the substitution of $\mathrm{\overline{\alpha}(x)}$ for $\mathrm{w}$ is free, is
\[\mathrm{FT. \;\;\;\;\forall \alpha_{B(\alpha)} \exists x \, R(\overline{\alpha}(x)) \rightarrow \exists n \forall \alpha_{B(\alpha)} \exists x \leq n \, R(\overline{\alpha}(x))},\]
where $\mathrm{B(\alpha) \equiv \forall x \, \alpha(x) \leq \beta(\overline{\alpha}(x))}$. For the {\em binary fan theorem}, which is no weaker over {\bf IRA}, $\mathrm{B(\alpha) \equiv \forall x \alpha(x) \leq 1}$.  Troelstra \cite{tr1974} proved that the full fan theorem  is conservative over Heyting arithmetic.

FT justifies a principle of {\em fan induction} with $\mathrm{R(w)}$ as basis and an arbitrary inductive predicate $\mathrm{A(w)}$.
For the binary fan the general form is
\begin{multline*}
\mathrm{\forall \alpha_{B(\alpha)} \exists x R(\overline{\alpha}(x)) ~\&~ \forall w_{B(w)} (R(w) \rightarrow A(w))} \\
\mathrm{ ~\&~ \forall w_{B(w)} (A(w \ast \langle 1 \rangle) ~\&~ A(w \ast \langle 2 \rangle) \rightarrow A(w)) \rightarrow A(1),}
\end{multline*}
where $\mathrm{B(\alpha) \equiv \forall x \, \alpha(x) \leq 1}$ and $\mathrm{B(w) \equiv \forall n < lh(w) \, (1 \leq (w)_n \leq 2)}$.
Modern reverse constructive mathematics establishes equivalences between restricted versions of FT and classically correct theorems of intuitionistic mathematics (e.g. \cite{ka2019}).

\subsection{Bar induction with a bar defined by a characteristic function}
Kleene's strongest restriction on the basis predicate $\mathrm{R(w)}$ leads to his weakest version
\begin{multline*}
\mathrm{BI_1. \; \; \;  \forall \alpha \exists x \rho(\overline{\alpha}(x)) = 0 ~\&~ \forall w (Seq(w) ~\&~ \rho(w) = 0 \rightarrow A(w))} \\
\mathrm{ ~\&~ \forall w (Seq(w) ~\&~ \forall s A(w \ast \langle s+1 \rangle) \rightarrow A(w)) \rightarrow A(1)}
\end{multline*}
($^x$26.3b in \cite{klve1965}) of bar induction.  Over {\bf IRA} this restriction is equivalent to requiring $\mathrm{R(w)}$ to be quantifier-free.  Solovay \cite{jrmrms} proved in primitive recursive arithmetic that  Kleene's {\bf I} + MP$_1$ is consistent relative to its subsystem {\bf IRA} + BI$_1$ + MP$_1$.

The corresponding version of the binary fan theorem is
\[\mathrm{FT_1. \;\;\;\;\forall \alpha_{B(\alpha)} \exists x \, \rho(\overline{\alpha}(x)) = 0 \rightarrow \exists n \forall \alpha_{B(\alpha)} \exists x \leq n \, \rho(\overline{\alpha}(x)) = 0},\]
where $\mathrm{B(\alpha) \equiv \forall x \, \alpha(x) \leq 1}$.  In Theorem 9.6 and Corollary 9.8 of \cite{vel2014}, Veldman has compiled a long list of theorems of intuitionistic mathematics equivalent to FT$_1$ over his minimal formal system {\bf BIM} (comparable to {\bf IRA}).  In particular, he showed that FT$_1$ is equivalent to the version of FT with $\mathrm{R(w) \equiv \exists n \, \beta(n) = w+1}$.

Kleene proved in \cite{klve1965} that the recursive sequences do not provide a classical omega-model of {\bf IRA} + FT$_1$ but the arithmetical sequences do; this distinction is exploited in \cite{vel2014}. BI$_1$ is stronger than FT$_1$ over {\bf IRA}; even the hyperarithmetical  sequences fail to satisfy BI$_1$.

\subsection{Decidable, thin and monotone bar induction}
Kleene formulated four axiom schemas ($^x$26.3a-d in \cite{klve1965}) of bar induction, including BI$_1$.

{\em Decidable bar induction} BI$_\mathrm{d}$ ($^x$26.3a) adds $\mathrm{\forall w (Seq(w) \rightarrow R(w) \vee \neg R(w))}$ to the hypotheses of BI. {\em Thin bar induction} BI! strengthens the assumption $\mathrm{\forall \alpha \exists x R(\overline{\alpha}(x))}$ of BI to  $\mathrm{\forall \alpha \exists ! x R(\overline{\alpha}(x))}$ for ($^x$26.3c), or to  $\mathrm{\forall \alpha \exists x (R(\overline{\alpha}(x)) ~\&~ \forall y<x \, \neg R(\overline{\alpha}(y)))}$ in the fourth version ($^x$26.3d).  BI$_\mathrm{d}$ is equivalent to BI! but stronger than BI$_1$ over {\bf IRA}.

Using BI! and continuous choice Kleene derived a fifth version, {\em monotone bar induction} BI$_\mathrm{mon}$ ($^\ast$27.13 in \cite{klve1965}), which adds
$\mathrm{\forall \alpha \forall x (R(\overline{\alpha}(x)) \rightarrow \forall y_{y>x} R(\overline{\alpha}(y)))}$ to the hypotheses of BI.
It was shown in \cite{jrmgvf2021} that BI$_\mathrm{d}$, BI$_\mathrm{mon}$ and BI$^\circ$ have the same classical content over {\bf IA$_1$}. From the classical point of view, BI$_\mathrm{d}$ and BI$_\mathrm{mon}$ express the full bar theorem (which is inconsistent with {\bf I}), but over {\bf IA$_1$} their negative interpretations are equivalent to (BI$^\circ$)$^g$ which is consistent with {\bf I}.

The corresponding versions FT$_\mathrm{d}$, FT! and FT$_\mathrm{mon}$ of the fan theorem are not all equivalent over {\bf IRA}.  Each version justifies a principle of restricted fan induction.  J. Berger (\cite{jbe2006}) proved that a special case c-FT of the monotone fan theorem is  constructively equivalent over {\bf HA$^\omega$} to the theorem that every pointwise continuous function from $\{0,1\}^\mathbb{N}$ to $\mathbb{N}$ is uniformly continuous.

\section{Two families of intuitionistically dubious principles}
If {\bf S} is a subsystem of {\bf B} then {\bf S$^\circ$} $\equiv$ {\bf S} + $\mathrm{(\neg \neg A \rightarrow A)}$ has the same language and mathematical axioms as {\bf S}, and {\bf S}$^{+g}$ $\subseteq$ {\bf S$^\circ$}; in this sense {\bf S$^\circ$} is to {\bf S} as {\bf PA} is to {\bf HA}. If it happens that {\bf S}$^{+g}$ $\not \subseteq$ {\bf S}, we seek an elegant characterization of the difference.

\subsection{Double negation shift principles}
\subsubsection{Double negation shift for numbers} This is the schema
\[\mathrm{DNS_0. \; \; \; \forall x \neg \neg A(x) \rightarrow \neg \neg \forall x A(x)}\]
for all formulas $\mathrm{A(x)}$ of the language. The converse is provable in {\bf IA$_1$}, so the $\mathrm{\rightarrow}$ can be strengthened to $\mathrm{\leftrightarrow}$.  {\bf IA$_1$} proves the restriction DNS$_0^-$ of DNS$_0$ to negative formulas $\mathrm{A(x)}$ since {\bf IA$_1$} $\vdash$ $\mathrm{\neg \neg A \leftrightarrow A}$ for every formula $\mathrm{A}$ not containing $\mathrm{\vee}$ or $\mathrm{\exists}$.

The restriction of DNS$_0$ to $\Sigma^0_1$ formulas $\mathrm{A(x)}$ is a weak consequence
\[\mathrm{{\Sigma^0_1}\mbox{-}DNS_0. \; \; \; \forall x \neg \neg \exists y \alpha(\langle x,y \rangle) = 0 \rightarrow \neg \neg \forall x \exists y \alpha(\langle x,y \rangle) = 0}\]
of MP$_1$ which Brouwer used in 1918 to prove that the intuitionistic real numbers form a closed species. Van Atten \cite{va2018} notes that Brouwer later formulated a stronger definition of ``closed'' in order to avoid this use of (a consequence of) Markov's Principle.

In \cite{scve1983} Scedrov and Vesley studied a principle of which $\Sigma^0_1$-DNS$_0$ is a special case.  They proved that {\bf B} $\not\vdash$ $\Sigma^0_1$-DNS$_0$ because $\Sigma^0_1$-DNS$_0$ fails in Krol's model of intuitionistic analysis \cite{krol1978}, and that {\bf B} + $\Sigma^0_1$-DNS$_0$ $\not\vdash$ MP$_1$. Their second argument, by $_S$realizability, establishes that {\bf I} + $\Sigma^0_1$-DNS$_0$ is consistent with Vesley's Schema VS which proves Brouwer's creating-subject counterexamples.

\subsubsection{Double negation shift and the negative interpretation of countable choice}

Over {\bf IA$_1$} + AC$_{01}$, (AC$_{01}$)$^g$ is equivalent to
\[\mathrm{({\Sigma^1_1 {neg}})\mbox{-}DNS_0. \; \; \; \forall x \neg \neg \exists \alpha R(x,\alpha) \rightarrow \neg \neg \forall x \exists \alpha R(x,\alpha)}\]
where $\mathrm{R(x,\alpha)}$ may be any negative formula, with parameters of both types allowed. The easy argument uses the fact that the converse of AC$_{01}$ is provable in {\bf IA$_1$}.

($\Sigma^1_1$neg)-DNS$_0$ is stronger than $\Sigma^0_1$-DNS$_0$, but still consistent with {\bf I} + VS because (AC$_{01}$)$^g$ is $_S$realizable (because it is negative and presumably true in classical Baire space). Hence ($\Sigma^1_1$neg)-DNS$_0$ is also $_S$realizable by Theorem 11.3(a) of \cite{klve1965}.

\subsubsection{Double negation shift for functions}
Full double negation shift for functions conflicts with Brouwer's continuity principles, but the version
\[\mathrm{DNS_1. \; \; \; \forall \alpha \neg \neg \exists x R(\overline{\alpha}(x)) \rightarrow \neg \neg \forall \alpha \exists x R(\overline{\alpha}(x))}\]
does not.  Lemma 27 in \cite{jrm2002} establishes that DNS$_1$, like MP$_1$, is self-realizing over {\bf B}, so {\bf I} + DNS$_1$ is consistent.  The useful special case
\[\mathrm{\Sigma^0_1\mbox{-}DNS_1. \; \; \; \forall \alpha \neg \neg \exists x \rho(\overline{\alpha}(x)) = 0 \rightarrow \neg \neg \forall \alpha \exists x \rho(\overline{\alpha}(x)) = 0}\]
is consistent with {\bf I} + VS by classical $_S$realizability, so might be considered in this context to be a palatable substitute for Markov's Principle.  Scedrov and Vesley observed in effect that {\bf IRA} + $\Sigma^0_1$-DNS$_1$ $\vdash$ $\Sigma^0_1$-DNS$_0$.

G\"odel, Dyson and Kreisel proved that the weak completeness of intuitionistic predicate logic for Beth semantics is equivalent, over {\bf IRA}, to a weaker consequence of $\Sigma^0_1$-DNS$_1$ which could be called the ``G\"odel-Dyson-Kreisel Principle'':\footnote{A doubly negated version
$\mathrm{\neg \neg WKL \equiv \forall n \exists \beta_{B(\beta)} \forall x \leq n \, \rho(\overline{\beta}(x)) \neq 0 \rightarrow \neg \neg \exists \beta_{B(\beta)} \forall x \rho(\overline{\beta}(x)) \neq 0}$ of weak K\"onig's Lemma is equivalent over {\bf IA$_1$} to $\neg \neg$FT$_1$ + GDK. A proof is in \cite{jrm????}, forthcoming.}
\[\mathrm{GDK. \; \; \; \forall \alpha_{B(\alpha)} \neg \neg \exists x \rho(\overline{\alpha}(x)) = 0 \rightarrow \neg \neg \forall \alpha_{B(\alpha)}\exists x \rho(\overline{\alpha}(x)) = 0}.\]
Because GDK is $^{\Delta^1_1}$realizable (\cite{jrm2010}) while $\Sigma^0_1$-DNS$_1$ is not, {\bf I} + GDK $\not\vdash$ $\Sigma^0_1$-DNS$_1$.
\subsection{Doubly negated characteristic function principles}
A number-theoretic relation $\mathrm{A(x)}$ (perhaps with number and sequence parameters) has a characteristic function for $\mathrm{x}$ only if it satisfies $\mathrm{\forall x (A(x) \vee \neg A(x))}$.
The doubly negated characteristic function (comprehension) schema
\[\mathrm{\neg \neg \, CF_0.  \; \; \; \neg \neg \exists \zeta \forall x (\zeta(x) = 0 \leftrightarrow A(x))}\]
says only that it is {\em persistently consistent} to assume a characteristic function for $\mathrm{A(x)}$ exists.  If {\bf S} proves an instance of $\mathrm{\neg\neg CF_0}$ in which the $\mathrm{A(x)}$ contains only $\mathrm{x}$ free, every consistent extension of {\bf S} is consistent with $\mathrm{\exists \zeta \forall x (\zeta(x) = 0 \leftrightarrow A(x))}$.

By Vafeiadou's characterization, the restriction $\mathrm{\neg \neg \, CF_0^{neg}}$ of $\mathrm{\neg \neg \, CF_0}$ to negative formulas $\mathrm{A(x)}$ is provable in the minimum classical extension of {\bf IA$_1$} + AC$_{00}!$.
An important special case\footnote{Over {\bf IRA} + CF$_\mathrm{d}$  or {\bf EL} + CF$_\mathrm{d}$, $\mathrm{\neg \neg \, \Pi^0_1\mbox{-}CF_0}$ is equivalent to the principle $\mathrm{\neg \neg \, \Pi^0_1\mbox{-}LEM}$ in \cite{fuko2018}, and $\mathrm{\neg \neg \Sigma^0_1}$-CF$_0$ is equivalent to $\mathrm{\neg \neg} \, \Sigma^0_1$-LEM.}, equivalent by intuitionistic logic to ($\Pi^0_1$-CF$_0$)$^g$, is
\[\mathrm{\neg \neg \Pi^0_1\mbox{-}CF_0. \; \; \; \; \forall \alpha \neg \neg \exists \zeta \forall x (\zeta(x) = 0 \leftrightarrow \forall y \alpha(\langle x,y \rangle) = 0)}.\]

\section{Minimum classical extensions of some subsystems of {\bf B}}

The negative translations of classical logical axioms and rules are correct by ituitionistic logic, so if E follows from $\mathrm{\Gamma}$ by classical logic then E$^g$ follows from $\mathrm{\Gamma}^g$ by intuitionistic logic.  With classical logic, E and E$^g$ are equivalent.  Even with intuitionistic logic, $\mathrm{\neg \neg E}^g$ and E$^g$ are equivalent.  These facts will be used without much comment in the following proofs.

\subsection{Theorem.}
\begin{enumerate}
\item[(i)]{{\bf (IA$_1$)$^{+g}$} = {\bf IA$_1$}.}
\item[(ii)]{{\bf (IRA)$^{+g}$} $\equiv$ ({\bf IA$_1$} + QF-AC$_{00}$)$^{+g}$ = {\bf IRA} + $\Sigma^0_1$-DNS$_0$.}
\item[(iii)]{({\bf IA$_1$} + AC$^{Ar}_{00}$)$^{+g}$ = {\bf IA$_1$} + AC$^{Ar}_{00}$ + $\Sigma^0_1$-DNS$_0$ + $\mathrm{\neg \neg \,\Pi^0_1\mbox{-}CF_0}$.}
\item[(iv)]{({\bf IA$_1$} + AC$_{00}!$)$^{+g}$ = {\bf IA$_1$} + AC$_{00}!$ + $\Sigma^0_1$-DNS$_0$ + $\mathrm{\neg \neg \, CF_0^\mathrm{neg}}$.}
\item[(v)]{({\bf IA$_1$} + AC$_{00}$)$^{+g}$ = {\bf IA$_1$} + AC$_{00}$ + $\Sigma^0_1$-DNS$_0$ + $\mathrm{\neg \neg \, CF_0^\mathrm{neg}}$.}
\item[(vi)]{({\bf IA$_1$} + AC$_{01}$)$^{+g}$ = {\bf IA$_1$} + AC$_{01}$ + $\mathrm{({\Sigma^1_1 {neg}})\mbox{-}DNS_0}$ = {\bf IA$_1$} + (AC$_{01}$)$^g$.}
\item[(vii)]{({\bf IA$_1$} + FT$_1$)$^{+g}$ = {\bf IA$_1$} + FT$_1$ + GDK.}
\item[(viii)]{({\bf IRA} + FT$_1$)$^{+g}$ = {\bf IRA} + FT$_1$ + $\Sigma^0_1$-DNS$_0$ + GDK.}
\item[(ix)]{({\bf IRA} + BI$_1$)$^{+g}$ = ({\bf IRA})$^{+g}$ + BI$_1$ + (BI$_1$)$^g$ $\subseteq$ {\bf IRA} + BI$_1$ + $\Sigma^0_1$-DNS$_1$.}
\item[(x)]{({\bf IA$_1$} + AC$_{00}$ + BI$_1$)$^{+g}$ = {\bf IA$_1$} + AC$_{00}$ + BI$_1$ + $\Sigma^0_1$-DNS$_0$ + $\mathrm{\neg \neg \, CF_0^\mathrm{neg}}$.}
\item[(xi)]{{\bf B}$^{+g}$ = ({\bf IA$_1$} + AC$_{01}$ + BI$_1$)$^{+g}$ = {\bf B} + $\mathrm{({\Sigma^1_1 {neg}})\mbox{-}DNS_0}$ = {\bf B} + (AC$_{01}$)$^g$.}
\end{enumerate}

{\em Proofs}. (i): The Gentzen negative translations of the axioms of {\bf IA$_1$} are provable in {\bf IA$_1$}, and the negative translations of the rules of inference are admissible for {\bf IA$_1$}, so no additions are needed.

(ii): To each quantifier-free formula $\mathrm{A(x,y)}$ there is by \cite{kl1969} a term $\mathrm{s(x,y)}$, with the same free variables, such that {\bf IA$_1$} proves both $\mathrm{\forall x \forall y (A(x,y) \leftrightarrow s(x,y) = 0)}$ and $\mathrm{\forall x \forall y (u(\langle x,y \rangle) = 0 \leftrightarrow s(x,y) = 0)}$ where $\mathrm{u = \lambda z . s((z)_0,(z)_1)}$.  Therefore {\bf IA$_1$} proves
$\mathrm{\exists \beta \forall x \forall y [A(x,y) \leftrightarrow \beta(\langle x,y \rangle) = 0]}$.  By intuitionistic logic the negative translation of $\mathrm{\forall x \exists y \beta(\langle x,y \rangle) = 0}$ is
equivalent to $\mathrm{\forall x \neg \neg \exists y \beta(\langle x,y \rangle) = 0}$, and the negative translation of
$\mathrm{\exists \alpha \forall x \beta(\langle x,\alpha(x) \rangle) = 0}$ is equivalent to
$\mathrm{\neg \neg \exists \alpha \forall x \beta(\langle x,\alpha(x) \rangle) = 0}$; therefore {\bf IRA} + $\Sigma^0_1$-DNS$_0$ $\vdash$ (QF-AC$_{00}$)$^g$. Conversely, $\Sigma^0_1$-DNS$_0$ is equivalent over {\bf IRA} to the negative translation of an instance of QF-AC$_{00}$.

(iii): Since QF-AC$_{00}$ is a special case of AC$^{Ar}_{00}$, {\bf IRA} $\subseteq$ {\bf IA$_1$} + AC$^{Ar}_{00}$.  By formula induction, {\bf IRA} + $\mathrm{\neg \neg \, \Pi^0_1\mbox{-}CF_0}$ proves $\mathrm{\neg \neg \exists \eta \forall x \forall y (\eta(\langle x,y \rangle) = 0 \leftrightarrow A(x,y))}$ for every negative arithmetical formula $\mathrm{A(x,y)}$.  The negative translation of AC$^{Ar}_{00}$ now follows using QF-AC$_{00}$ and $\Sigma^0_1$-DNS$_0$ as in (ii).  This is a variation of Solovay's argument;  he started with MP$_1$ and
$\mathrm{\neg \neg \, \Sigma^0_1\mbox{-}CF_0}$ instead of $\mathrm{\Sigma^0_1\mbox{-}DNS_0}$ and $\mathrm{\neg \neg \, \Pi^0_1\mbox{-}CF_0}$, which give a precise characterization here.  See the next theorem also.

Conversely, $\mathrm{\forall x \exists z (z = 0 \leftrightarrow \forall y \alpha(\langle x,y \rangle) = 0) \rightarrow
\exists \zeta \forall x (\zeta(x) = 0 \leftrightarrow \forall y \alpha(\langle x,y \rangle) = 0)}$ is an instance of AC$^{Ar}_{00}$, and
$\mathrm{\forall x \neg \neg \exists z (z = 0 \leftrightarrow \forall y \alpha(\langle x,y \rangle) = 0)}$ is provable in {\bf IA$_1$}.  It follows that $\mathrm{\neg \neg \exists \zeta \forall x (\zeta(x) = 0 \leftrightarrow \forall y \alpha(\langle x,y \rangle) = 0)}$ is provable in {\bf IA$_1$} + (AC$^{Ar}_{00}$)$^g$.

(iv): {\bf IA$_1$} + AC$_{00}!$ = {\bf IRA} + CF$_{\mathrm{d}}$ by Vafeiadou's characterization; therefore ({\bf IA$_1$} + AC$_{00}!$)$^{+g}$ =
({\bf IRA} + CF$_{\mathrm{d}}$)$^{+g}$ = ({\bf IRA})$^{+g}$ + CF$_{\mathrm{d}}$ + (CF$_{\mathrm{d}}$)$^g$.  Each instance of $\mathrm{\neg \neg \, CF^\mathrm{neg}_0}$ is equivalent over {\bf IA$_1$} to the conclusion of the negative translation of an instance of CF$_{\mathrm{d}}$, and the negative translation $\mathrm{\forall x \neg (\neg A}^g\mathrm{(x) ~\&~ \neg \neg A}^g\mathrm{(x))}$
of $\mathrm{\forall x (A(x) \vee \neg A(x))}$ is provable in {\bf IA$_1$} for all formulas $\mathrm{A(x)}$, so  (CF$_{\mathrm{d}}$)$^g$ and $\mathrm{\neg \neg \, CF^\mathrm{neg}_0}$ are equivalent over {\bf IA$_1$}.

(v) follows from (iv) because AC$_{00}$ and AC$_{00}$! are equivalent as schemas over {\bf IA$_1$} + $\mathrm{(\neg \neg A \rightarrow A)}$, which proves $\mathrm{\forall x (\exists y A(x,y) \rightarrow \exists ! y (A(x,y) ~\&~ \forall z < y \neg A(x,z))}$. Therefore (AC$_{00}$)$^g$ and (AC$_{00}!$)$^g$ are equivalent over {\bf IA$_1$}.

(vi) is immediate from the definitions.

(vii): It is routine to show that {\bf IA$_1$} + FT$_1$ + GDK proves (FT$_1$)$^g$. The proof of GDK in {\bf IA$_1$} + (FT$_1$)$^g$ is an easy  exercise.  (viii) follows by (ii).

(ix):  It is routine to show that {\bf IRA} + BI$_1$ + $\Sigma^0_1$-DNS$_1$ proves (BI$_1$)$^g$, and $\Sigma^0_1$-DNS$_0$ follows from $\Sigma^0_1$-DNS$_1$ in {\bf IRA}.  Now use (ii).

(x) and (xi) follow from (v) and (vi) because {\bf IA$_1$} + AC$_{00}$ + $\mathrm{(\neg \neg A \rightarrow A)}$ $\vdash$ BI$_1$ (cf. $^\ast$26.1$^\circ$ in \cite{klve1965}), so ({\bf IA$_1$} + AC$_{00}$)$^{+g}$ $\vdash$ (BI$_1$)$^g$. \qed

\subsection{Corollary.} {\bf IRA} + BI$_1$ + $\Sigma^0_1$-DNS$_1$ is its own minimum classical extension.

\vskip 0.1cm

{\em Proof.} By Theorem 5.1(ix) with the observation that {\bf IA$_1$} $\vdash$ ($\Sigma^0_1$-DNS$_1$)$^g$. \qed
\subsection{Corollary.}For each subsystem {\bf S} of {\bf B} considered in Theorem 5.1:
\begin{enumerate}
\item[(i)]{{\bf S}$^{+g}$ is its own minimum classical extension.}
\item[(ii)]{({\bf S} + MP$_1$)$^{+g}$ = {\bf S}$^{+g}$ + MP$_1$ is its own minimum classical extension.}
\item[(iii)]{{\bf S}$^{+g}$ + MP$_1$ is consistent with strong continuous choice CC$_{11}$ ($^{\mathrm{x}}$27.1 in \cite{klve1965}).}
\item[(iv)]{{\bf S}$^{+g}$ + CC$_{11}$ $\not\vdash$ MP$_1$.}
\end{enumerate}

{\em Proofs.} (i) is true because the Gentzen negative translation is idempotent.  (ii) is true because {\bf S} $\vdash$ (MP$_1$)$^g$.  The rest is implicit in \cite{klve1965}. (ii) holds by classical Kleene function-realizability (cf. Lemma 8.4(a) of \cite{klve1965}).  (iv) holds because every theorem of {\bf S}$^{+g}$ + CC$_{11}$ is $_S$realizable but MP$_1$ is not (cf. Lemma 10.7, Theorem 11.3 and Corollary 11.10(a) in \cite{klve1965}). \qed

\subsection{Corollary} Each of {\bf IRA} + MP$_1$, {\bf IA$_1$} + FT$_1$ + MP$_1$, {\bf IRA} + FT$_1$ + MP$_1$ and {\bf IRA} + BI$_1$ + MP$_1$ is its own minimum classical extension.

{\em Proof.} {\bf IA$_1$} + MP$_1$ proves  $\Sigma^0_1$-DNS$_0$, $\Sigma^0_1$-DNS$_1$ and GDK so the results follows from Theorem 5.1(ii), (vii), (viii) and (ix) using Corollary 5.3.

\subsection{Two questions} Sometimes only one or two additional axioms must be added to a subsystem {\bf S} of {\bf B} in order to prove its G\"odel-Gentzen negative interpretation.  The unrestricted axioms of countable choice and comprehension have resisted this treatment, requiring instead the addition of an axiom schema $\mathrm{\neg \neg \, CF_0^\mathrm{neg}}$ or $\mathrm{({\Sigma^1_1 {neg}})\mbox{-}DNS_0}$. Is there a more elegant solution?

$\Sigma^0_1$-DNS$_1$ evidently suffices for the negative interpretation of BI$_1$, but is it stronger than necessary?  Does {\bf IRA} + BI$_1$ + $\Sigma^0_1$-DNS$_0$ + (BI$_1$)$^g$ $\vdash$ $\Sigma^0_1$-DNS$_1$?

\section{Bar induction in two contexts}

The next result sharpens Solovay's proof that {\bf IA$_1$} + AC$^{Ar}_{00}$ + BI$_1$ + $\mathrm{(\neg \neg A \rightarrow A)}$ can be negatively interpreted in {\bf IRA} + BI$_1$ + MP$_1$.  In fact he proved the stronger theorem (cf. \cite{jrmrms}) that $\mathrm{\neg \neg \, \Sigma^0_1\mbox{-}CF_0}$ (thus  $\mathrm{\neg \neg \, CF_0^{Ar}}$) holds in {\bf IRA} + BI$_1$ + MP$_1$, but $\mathrm{\neg \neg \, \Pi^0_1\mbox{-}CF_0}$ gives the double negation of the characteristic function principle for {\em negative} arithmetical formulas, which suffices with $\Sigma^0_1$-DNS$_0$ and QF-AC$_{00}$ for the negative interpretation of AC$^{Ar}_{00}$.  For the derivation of $\mathrm{\neg \neg \, \Pi^0_1\mbox{-}CF_0}$ by bar induction and for the negative interpretation of BI$_1$, MP$_1$ is not needed; $\Sigma^0_1$-DNS$_1$ suffices.

\subsection{Theorem.} (after Solovay)
\begin{enumerate}
\item[(i)]{{\bf IA$_1$} + (BI$_1$)$^g$  $\vdash$ $\mathrm{\neg \neg \, \Pi^0_1\mbox{-}CF_0}$.}
\item[(ii)]{{\bf IA$_1$} +  AC$^{Ar}_{00}$ + BI$_1$ + $\mathrm{(\neg \neg A \rightarrow A)}$ can be negatively interpreted in (and therefore is equiconsistent with) its subsystem {\bf IRA} + BI$_1$ + $\Sigma^0_1$-DNS$_1$.}
\end{enumerate}

{\em Proofs.}
(i): Adapting Solovay's argument that {\bf IA$_1$} + BI$_1$ + MP$_1$ $\vdash$ $\mathrm{\neg \neg \, \Sigma^0_1\mbox{-}CF_0}$ (as in \cite{jrm2003}, \cite{jrmrms}), assume for contradiction {\bf (a)} $\mathrm{\forall \zeta \neg \forall x (\zeta(x) = 0 \leftrightarrow \forall y \alpha(\langle x,y \rangle) = 0)}$.  Then {\bf (b)} $\mathrm{\forall \zeta \neg \neg \exists x [(\zeta(x) = 0 ~\&~ \neg \neg \exists y \alpha(\langle x,y \rangle) \neq 0) \vee (\zeta(x) \neq 0 ~\&~ \forall y \alpha(\langle x,y \rangle) = 0)]}$ follows in {\bf IA$_1$}, and this entails
{\bf (c)} $\mathrm{\forall \zeta \neg \neg \exists x [(\zeta((x)_0) = 0 ~\&~ \alpha(\langle (x)_0,(x)_1 \rangle) \neq 0) ~\vee}$

\noindent $\mathrm{(\zeta((x)_0) \neq 0 ~\&~ \forall y \alpha(\langle (x)_0,y \rangle) = 0)]}$.

In {\bf IRA} one can define a binary sequence $\mathrm{\rho}$ such that $\mathrm{\rho(w) = 0}$ if and only if $\mathrm{Seq(w)}$ and for some $\mathrm{j < lh(w)}$ either
\begin{enumerate}
\item[{\bf (d)}]{$\mathrm{(w)_j = 1 ~\&~ \exists y < lh(w)\, \alpha(\langle j,y \rangle) \neq 0}$, or}
\item[{\bf (e)}]{$\mathrm{(w)_j > 1 ~\&~ [\alpha(\langle j, ((w)_j\stackrel{.}{-}2) \rangle) = 0 \vee \exists y < (w)_j\stackrel{.}{-}2 \; \alpha(\langle j,y \rangle) \neq 0]}$.}
\end{enumerate}
Now prove {\bf (f)} $\mathrm{\forall \zeta \neg \neg \exists n \rho(\overline{\zeta}(n)) = 0}$ by cases on (c) using (d) and (e), giving the first hypothesis for an application of (BI$_1$)$^g$. The negative inductive predicate $\mathrm{A(w)}$ is
\begin{multline*}
\mathrm{A(w) \equiv \neg \neg \exists j < lh(w) [((w)_j = 1 \rightarrow \neg \forall y \alpha(\langle j,y \rangle) = 0)} \\ \mathrm{\&~ ((w)_j > 1 \rightarrow [\alpha(\langle j, ((w)_j\stackrel{.}{-}2) \rangle) \neq 0 \rightarrow \exists y < ((w)_j\stackrel{.}{-}2) \; \alpha(\langle j,y \rangle) \neq 0])].}
\end{multline*}
Evidently {\bf (g)} $\mathrm{\forall w (Seq(w) ~\&~ \rho(w) = 0 \rightarrow A(w))}$.  In order to establish the inductive hypothesis
{\bf (h)} $\mathrm{\forall w (Seq(w) ~\&~ \forall s A(w*\langle s+1\rangle) \rightarrow A(w))}$, argue by contradiction as follows, noting that
in general $\mathrm{(w*\langle n \rangle)_{lh(w)} = n}$.

Assume $\mathrm{Seq(w) ~\&~ \forall s A(w*\langle s+1\rangle) ~\&~ \neg A(w)}$.  From $\mathrm{A(w*\langle 1 \rangle)}$ and $\mathrm{\neg A(w)}$ we get $\mathrm{(w*\langle 1 \rangle)_{lh(w)} = 1 ~\&~ \neg \forall y \alpha(\langle lh(w),y \rangle) = 0}$.  From $\mathrm{\forall n A(w*\langle n+2 \rangle)}$ and $\mathrm{\neg A(w)}$  we get $\mathrm{\forall n [(w*\langle n+2 \rangle)_{lh(w)} > 1 ~\&~ (\alpha(\langle lh(w),n \rangle) \neq 0 \rightarrow \exists y < n \, \alpha(\langle lh(w),y \rangle) \neq 0)}]$, from which it follows that
$\mathrm{\forall n (\alpha(\langle lh(w),n \rangle) \neq 0 \rightarrow \exists y < n \, \alpha(\langle lh(w),y \rangle) \neq 0)}$, contradicting
$\mathrm{\neg \forall y \alpha (\langle lh(w),y \rangle) = 0}$.  This completes the proof of (h).

By (BI$_1$)$^g$ conclude  $\mathrm{A(\langle \, \rangle)}$, which is impossible because $\mathrm{lh(\langle \, \rangle) = 0}$.  Therefore $\mathrm{\neg \neg \, \Pi^0_1\mbox{-}CF_0}$ holds in {\bf IRA} + (BI$_1$)$^g$.

(ii): (BI$_1$)$^g$ was treated in Theorem 5.1(ix), and {\bf IRA} + (BI$_1$)$^g$ $\vdash$ (AC$^{Ar}_{00}$)$^g$ by formula induction from (i) (cf. the proof of Theorem 5.1(iii)).  Observe that {\bf IRA} $\subseteq$ {\bf IA$_1$} + AC$_{00}^{Ar}$, and {\bf IA$_1$} proves $\mathrm{(\neg \neg A \rightarrow A)}^g$ for all formulas $\mathrm{A}$. \qed

\subsection{Theorem.} {\bf IRA} + (BI$_1$)$^g$ + $\mathrm{\neg \neg \, \Pi^1_1\mbox{-}CF_0}$ $\vdash$ $\Sigma^0_1$-DNS$_1$, where
$\mathrm{\neg \neg \Pi^1_1\mbox{-}CF_0}$ is
\[\mathrm{\forall \gamma \neg \neg \exists \zeta \forall x (\zeta(x) = 0 \leftrightarrow \forall \alpha \exists y \, \gamma(\overline{\alpha}(\langle x,y \rangle)) = 0).}\]

{\em Proof.} Assume {\bf (a)} $\mathrm{\forall \alpha \neg \neg \exists x \rho(\overline{\alpha}(x)) = 0}$.  The goal is to prove $\mathrm{\neg \neg \forall \alpha \exists x \rho(\overline{\alpha}(x)) = 0}$ in {\bf IA$_1$} + $\mathrm{\neg \neg \, \Pi^1_1\mbox{-}CF_0}$ + (BI$_1$)$^g$.   First define in {\bf IA$_1$} a function $\mathrm{\gamma}$ such that {\bf (b)} $\mathrm{\forall \alpha \forall x (Seq(x) \rightarrow \forall y(\rho(x * \overline{\alpha}(y)) = 0 \leftrightarrow \gamma(\overline{\alpha}(\langle x,y \rangle)) = 0))}$.
As  the desired conclusion is negative, assume for ``$\mathrm{\neg \neg \exists}$-elimination'' (cf. \cite{jrm2003}) from the appropriate instance of $\mathrm{\neg \neg \, \Pi^1_1\mbox{-}CF_0}$: {\bf (c)} $\mathrm{\forall x (\zeta(x) = 0 \leftrightarrow \forall \alpha \exists y \gamma(\overline{\alpha}(\langle x,y \rangle)) = 0)}$.  Then in particular {\bf (d)} $\mathrm{\forall w (Seq(w) \rightarrow (\zeta(w) = 0 \leftrightarrow \forall \alpha \exists y \rho(w * \overline{\alpha}(y)) = 0))}$.

From (d) follow the other hypotheses {\bf(e)} $\mathrm{\forall w (Seq(w) ~\&~ \rho(w) = 0 \rightarrow \zeta(w) = 0)}$ and {\bf (f)} $\mathrm{\forall w (Seq(w) ~\&~ \forall s \, \zeta(w * \langle s+1 \rangle) = 0 \rightarrow \zeta(w) = 0)}$ of the instance of (BI$_1$)$^g$ with $\mathrm{\zeta(w) = 0}$ as the inductive predicate, so {\bf (g)} $\mathrm{\zeta(\langle \, \rangle) = 0}$, hence $\mathrm{\forall \alpha \exists x \rho(\overline{\alpha}(x)) = 0}$.  Discharging hypothesis (c) by $\neg \neg \exists$-elimination,
{\bf (h)} $\mathrm{\neg \neg \forall \alpha \exists x \rho(\overline{\alpha}(x)) = 0}$.  \qed

\section{Minimum classical extensions of systems between {\bf B} and {\bf I}}
The principle of monotone bar induction BI$_\mathrm{mon}$ is provable in {\bf I} and in {\bf B$^\circ$} but not in {\bf B}, and {\bf IA$_1$} + BI$_\mathrm{mon}$ proves BI$_\mathrm{d}$.  It follows that the variant {\bf B$'$} of {\bf B} with BI$_\mathrm{mon}$ as an axiom schema in place of BI$_\mathrm{d}$ is classically sound and lies strictly between {\bf B} and {\bf I}.\footnote{Veldman's careful analysis of Brouwer's writing on the subject led him to the conclusion that Brouwer sometimes assumed a monotone bar, but sometimes fell into the error of trying to justify classical bar induction (which is inconsistent with his own continuity principle).} As it happens, {\bf B$'$} has the same classical content as {\bf B} over {\bf IA$_1$} (cf. \cite{jrmgvf2021}), but this may not be the case for every classically sound intermediate system.

\subsection{Neighborhood function principles}
A classically sound choice principle guaranteeing that every pointwise continuous relation has a modulus of continuity is Troelstra's {\em neighborhood function principle}:
\[\mathrm{NFP. \; \forall \alpha \exists x A(\overline{\alpha}(x)) \rightarrow \exists \sigma \forall \alpha [\exists ! x \sigma(\overline{\alpha}(x)) > 0 ~\&~ \forall x \forall y (\sigma(\overline{\alpha}(x)) = y+1 \rightarrow A(\overline{\alpha}(y)))]}.\]
This version, labeled AC$_{1/2,0}$ in \cite{jrmgvf2012}, is equivalent to Troelstra's over {\bf IA$_1$}.

NFP  follows easily from  ``Brouwer's Principle for a Number'' CC$_{10}$ ($^\ast$27.2 in \cite{klve1965}) but is not provable in {\bf B}.
The monotone version NFP$_\mathrm{mon}$ (AC$^\mathrm{m}_{1/2,0}$ in \cite{jrmgvf2012}) of NFP is interderivable with BI$_\mathrm{mon}$ over {\bf B} (so does not add classical content to {\bf B}), but NFP is apparently stronger.  A partial characterization of ({\bf IRA} + NFP)$^{+g}$ follows.

\subsection{Theorem}
\begin{enumerate}
\item[(i)]{({\bf IRA} + NFP)$^{+g}$ $\subseteq$ {\bf IRA} + NFP + $\Sigma^0_1$-DNS$_1$ + $\neg \neg \, $CF$_0^\mathrm{neg}$.}
\item[(ii)]{{\bf IRA} + NFP$^g$ + $\Sigma^0_1$-DNS$_0$ $\vdash$ $\neg \neg \, $CF$_0^\mathrm{neg}$.}
\end{enumerate}

{\em Proofs.} (i) Assume {\bf (a)} $\mathrm{\forall \alpha \neg \neg \exists x A}^g{(\overline{\alpha}(x))}$. For $\mathrm{\neg \neg \exists \zeta}$-elimination from $\neg\neg$CF$_0^\mathrm{neg}$: {\bf (b)} $\mathrm{\forall w [\zeta(w) = 0 \leftrightarrow A}^g\mathrm{(w)]}$. From (a), (b) by
$\Sigma^0_1$-DNS$_1$: {\bf (c)} $\mathrm{\neg \neg \forall \alpha \exists x \zeta(\overline{\alpha}(x)) = 0.}$
NFP gives {\bf (d)} $\mathrm{\neg \neg \exists \sigma \forall \alpha [\exists ! x \sigma(\overline{\alpha}(x)) > 0 ~\&~ \forall x \forall y (\sigma(\overline{\alpha}(x)) = y + 1 \rightarrow \zeta(\overline{\alpha}(y)) = 0)]}$, whence {\bf (e)} $\mathrm{\neg \neg \exists \sigma \forall \alpha [\exists ! x \sigma(\overline{\alpha}(x)) > 0 ~\&~ \forall x \forall y (\sigma(\overline{\alpha}(x)) = y + 1 \rightarrow A}^g\mathrm{(\overline{\alpha}(y))]}$ by (b) and {\em a fortiori} {\bf (f)} $\mathrm{\neg \neg \exists \sigma \forall \alpha [\neg \neg \exists ! x \sigma(\overline{\alpha}(x)) > 0 ~\&~ \forall x \forall y (\sigma(\overline{\alpha}(x)) = y + 1 \rightarrow A}^g\mathrm{(\overline{\alpha}(y))]}$. Because (f) is a negation not involving $\mathrm{\zeta}$, (b) may now be discharged by $\mathrm{\neg \neg \exists \zeta}$-elimination.

(ii) {\bf IA$_1$} $\vdash$ $\mathrm{\forall x \neg \neg (A(x) \vee \neg A(x))}$ and so {\bf (a)}
$\mathrm{\forall x \neg \neg \exists y (y \leq 1 ~\&~ (y = 0 \leftrightarrow A(x)))}$.  Let B(w) abbreviate
$\mathrm{Seq(w) ~\&~ 1 \leq lh(w) \leq 2 ~\&~ (lh(w) = 1 \leftrightarrow A((w)_0 \stackrel{.}{-} 1))}$, where for $\neg \neg \, $CF$_0^\mathrm{neg}$ the
A(w) and therefore B(w) are negative. Then {\bf(b)} $\mathrm{\forall \alpha \neg \neg \exists y B(\overline{\alpha}(y))}$, so by NFP$^g$:
{\bf (c)} $\mathrm{\neg \neg \exists \sigma \forall \alpha [\neg \neg \exists ! x \sigma(\overline{\alpha}(x)) > 0 ~\&~ \forall x \forall y (\sigma(\overline{\alpha}(x)) = y+1 \rightarrow B(\overline{\alpha}(y)))]}$.  Now assume {\bf (d)} $\mathrm{\forall \alpha \neg \neg \exists ! x \sigma(\overline{\alpha}(x)) > 0 ~\&~ \forall \alpha \forall x \forall y (\sigma(\overline{\alpha}(x)) = y+1 \rightarrow B(\overline{\alpha}(y)))}$ for $\mathrm{\neg \neg \exists \sigma}$-elimination from (c), since $\neg \neg \, $CF$_0^\mathrm{neg}$ is negative. Substituting $\mathrm{\lambda t.n}$ for $\mathrm{\alpha}$ in (d) gives {\bf (e)} $\mathrm{\forall n \neg \neg \exists ! x \sigma(\overline{\lambda t.n}(x)) > 0 ~\&~ \forall n \forall x \forall y (\sigma(\overline{\lambda t. n}(x)) = y+1 \rightarrow B(\overline{\lambda t.n}(y)))}$, hence
{\bf(f)} $\mathrm{\neg \neg \forall n \exists x \sigma(\overline{\lambda t.n}(x)) > 0 ~\&~ \forall n \forall x \forall y (\sigma(\overline{\lambda t. n}(x)) = y+1 \rightarrow B(\overline{\lambda t.n}(y)))}$ by $\Sigma^0_1$-DNS$_0$, so by QF-AC$_{00}$: {\bf (g)} $\mathrm{\neg \neg \exists \tau \forall n \sigma(\overline{\lambda t.n}(\tau(n))) > 0}$.  For $\mathrm{\neg \neg \exists \tau}$-elimination from (g) assume
{\bf (h)} $\mathrm{\forall n \, \sigma(\overline{\lambda t.n}(\tau(n))) > 0}$, so {\bf (i)} $\mathrm{\forall n B(\overline{\lambda t. n}(\sigma(\overline{\lambda t.n}(\tau(n)))\stackrel{.}{-} 1))}$ by (f). It follows that {\bf (j)} $\mathrm{\forall n [1 \leq lh(\overline{\lambda t. n}(\sigma(\overline{\lambda t.n}(\tau(n)))\stackrel{.}{-} 1)) = \sigma(\overline{\lambda t.n}(\tau(n)))\stackrel{.}{-} 1 \leq 2]}$, so
{\bf (k)} $\mathrm{\forall n \,[(\overline{\lambda t. n}(\sigma(\overline{\lambda t.n}(\tau(n)))\stackrel{.}{-} 1)_0 \stackrel{.}{-} 1 = n]}$.
Finally set $\mathrm{\zeta = \lambda n.\sigma(\overline{\lambda t.n}(\tau(n)))\stackrel{.}{-} 2}$ to conclude $\mathrm{\exists \zeta \forall n (\zeta(n) = 0 \leftrightarrow A(n))}$.  Two $\mathrm{\neg \neg \exists}$-eliminations, discharging (h) and (d) respectively, complete the proof of $\neg \neg \, $CF$_0^\mathrm{neg}$. \qed

\subsection{Dependent choice for sequences}  Dependent choice for numbers DC$_0$ is a theorem of {\bf IA$_1$} + AC$_{00}$, but dependent choice for sequences
\[\mathrm{DC_1. \; \; \forall \alpha \exists \beta A(\alpha,\beta) \rightarrow \forall \alpha \exists \gamma [(\gamma)_0 = \alpha ~\&~ \forall n A((\gamma)_n,(\gamma)_{n+1})]}\]
(where $\mathrm{(\gamma)_n = \lambda x. \gamma(\langle n,x \rangle)}$)
is not obviously provable in {\bf B} or even in {\bf B$^\circ$}. It is not hard to see, however, that {\bf B} + DC$_1$ is a classically sound subsystem of {\bf I}.

\subsection{Theorem}
\begin{enumerate}
\item[(i)]{{\bf I} $\vdash$ DC$_1$.}
\item[(ii)]{{\bf IRA} + DC$_1$ $\vdash$ AC$_{01}$.}
\item[(iii)]{({\bf B} + DC$_1$)$^{+g}$ = {\bf IRA} + BI$_1$ + DC$_1$ + (DC$_1$)$^g$.}
\end{enumerate}

{\em Proofs.} (i) Assume {\bf (a)} $\mathrm{\forall \alpha \exists \beta A(\alpha,\beta)}$.  By CC$_{11}$ there is a $\mathrm{\sigma}$ satisfying
{\bf (b)} $\mathrm{\forall \alpha \exists ! \beta [\{\sigma\}[\alpha] \simeq \beta ~\&~ A(\alpha,\beta)]}$.  Fix $\mathrm{\alpha}$.  It will be enough to show that there is a $\mathrm{\zeta}$ such that
{\bf (c)} $\mathrm{\forall n [((\zeta)_n)_0 = \alpha ~\&~ \forall i < n [A(((\zeta)_n)_i,((\zeta)_n)_{i+1}) ~\&~ ((\zeta)_n)_i = ((\zeta)_{n+1})_i]]}$, since then we can define $\mathrm{\gamma}$ so that {\bf (d)} $\mathrm{\forall n [(\gamma)_n = ((\zeta)_{n+2})_n]}$ and then it will follow that $\mathrm{(\gamma)_0 = \alpha}$ and for all $\mathrm{n}$: $\mathrm{(\gamma)_{n+1} = ((\zeta)_{n+3})_{n+1} = ((\zeta)_{n+2})_{n+1}}$ so $\mathrm{A((\gamma)_n,(\gamma)_{n+1})}$.

Toward (c), first prove by induction:
$\mathrm{\forall n \exists \delta [(\delta)_0 = \alpha ~\&~ \forall i < n (\delta)_{i+1} \simeq \{\sigma\}[(\delta)_i]]}$. By AC$_{01}$,
$\mathrm{\exists \zeta \forall n [((\zeta)_n)_0 = \alpha ~\&~ \forall i < n ((\zeta)_n)_{i+1} \simeq \{\sigma\}[((\zeta)_n)_i]]}$. Apply (b).\footnote{The logic of partial terms is {\em not} involved in this argument because the informal expression $\mathrm{\{\sigma\}[\alpha]}$, which helps to clarify the proof, always designates a fully defined sequence $\mathrm{\beta}$ satisfying
$\mathrm{\forall x \forall y [\beta(x) = y \leftrightarrow \exists z [\sigma(\langle x+1 \rangle \ast \overline{\alpha}(z)) = y+1 ~\&~ \forall n<z \, \sigma(\langle x+1 \rangle \ast \overline{\alpha}(n)) = 0]]}$.}

(ii)  Assume {\bf (a)} $\mathrm{\forall n \exists \alpha A(n,\alpha)}$.  We want to show $\mathrm{\exists \beta \forall n A(n,(\beta)_n)}$.  From (a) we conclude {\bf (b)} $\mathrm{\forall \alpha \exists \beta [\beta(0) = \alpha(0)+1 ~\&~ A(\alpha(0),\lambda x.\beta(x+1))]}$, from which DC$_1$ gives {\bf (c)} $\mathrm{\exists \gamma [(\gamma)_0 = \lambda t. 0 ~\&~ \forall n [(\gamma)_{n+1}(0) = (\gamma)_n(0)+1 ~\&~ A((\gamma)_n(0),\lambda x. (\gamma)_{n+1}(x+1))]]}$.  For any such $\mathrm{\gamma}$, {\bf (d)} $\mathrm{\forall n \,(\gamma)_n(0) = n}$ and {\bf (e)} $\mathrm{\forall n \,  A((\gamma)_n(0),\lambda x. (\gamma)_{n+1}(x+1))}$ hold by induction, and it is easy to define a $\mathrm{\beta}$ such that $\mathrm{\forall n \forall x \, \beta(\langle n,x \rangle )= (\gamma)_{n+1}(x+1)}$.

(iii) is immediate from (ii) and (the proof of) Theorem 5.1(xi). \qed

\subsection{Bar induction of type one}
One year after the publication of \cite{klve1965}, Howard and Kreisel \cite{hokr1966} corrected a typo in the conclusion of the formulation, in Section 6.3 of  \cite{sp1962}, of Spector's axiom schema of bar induction of type one; call the corrected version ``SBI$_1$.'' In a minor variant {\bf H} of {\bf IA$_1$}, using the continuous choice principle CC$_{10}$ which is a theorem of {\bf I}, they derived SBI$_1$ from the special case  BI$_\mathrm{QF}$ of BI in which $\mathrm{R(w)}$ is required to be quantifier-free (hence decidable, even if free sequence variables are present).  Their proof also shows that {\bf IRA} + BI$_1$ + CC$_{10}$ $\vdash$ SBI$_1$.

In Appendix 1 of \cite{hokr1966} Howard and Kreisel proved that DC$_1$ is equivalent to SBI$_1$ over {\bf H$^\circ$} $\equiv$ {\bf H} + $\mathrm{(\neg \neg A \rightarrow A)}$, and they asked (Problem 9 of the appendix) whether DC$_1$ is derivable from AC$_{01}$ over {\bf H} or {\bf H$^\circ$}. The same questions can be asked over {\bf IRA} and {\bf IRA$^\circ$}.  Is it possible that {\bf B$^{+g}$} = ({\bf B} + SBI$_1$)$^{+g}$?

\subsection{A sticky question}
Is there a reasonable way to define a unique minimum classical extension of a consistent theory which is not classically consistent? How can one make sense of cls({\bf I}), and hence of {\bf I}$^{+g}$?

Gandy, Kreisel and Tait \cite{gakrta1960} proved that every sequence belonging to all classical $\omega$-models of {\bf B} is  hyperarithmetic, and therefore $\Delta^1_1$ by the Suslin-Kleene Theorem.  Kleene (\cite{kl1955b}, XXIV) proved that the $\Delta^1_1$ sequences do not satisfy BI$_1$.  Classical $\omega$-models of subsystems of {\bf I} extending {\bf B} must include all projective sequences.

The negative interpretations of principles like NFP, DC$_1$ and BI$_\mathrm{mon}$ which are provable in {\bf I}, and whose addition to {\bf B} is classically sound, should obviously belong to {\bf I}$^{+g}$.  But {\bf I} refutes very simple consequences of the law of excluded middle to which Bishop constructivists have given colorful names.  For example, {\bf I} proves ($^\ast$27.17 in \cite{klve1965}) the {\em negation} of the ``weak limited principle of omniscience''
\[\mathrm{WLPO. \; \; \; \forall \alpha (\forall x \alpha(x) = 0 \vee \neg \forall x \alpha (x) = 0).}\]

One possibility is suggested by the fact that {\bf I} is consistent with the collection of all {\em negative} sentences of the language $\mathcal{L}({\bf I})$ which are {\em true in classical Baire space} $\mathcal{B}$ = $(\omega, \omega^\omega)$.  Let $\mathcal{X}_\mathcal{B}$ be the collection of all subsystems {\bf S} of {\bf I} which extend {\bf B} and prove {\em only statements true in} $\mathcal{B}$. Then $\bigcup \{\mathrm{\bf S} : \mathrm{\bf S} \in \mathcal{X}_\mathcal{B}\}$ is classically sound, and proves just the theorems of {\bf I} which are true in $\mathcal{B}$.  The {\em minimum classical extension} {\bf I$_\mathcal{B}^{+g}$} {\em of} {\bf I} {\em relative to} $\mathcal{B}$ may be identified with  $\bigcup \{\mathrm{\bf S}^{+g} : \mathrm{\bf S} \in \mathcal{X}_\mathcal{B}\}$ + CC$_{11}$.

\subsubsection{\bf Theorem} (Vafeiadou, \cite{jrmgvf2021})  By this definition, {\bf I$_\mathcal{B}^{+g}$} = {\bf I} + $\mathrm{(\Gamma_\mathcal{B}^\circ)^g}$ where $\mathrm{\Gamma_\mathcal{B}^\circ}$ is the collection of all sentences in $\mathcal{L}({\bf I})$ which are true in $\mathcal{B}$.

\vskip 0.2cm

This appeal to truth in $\mathcal{B}$ appears necessary. If {\bf S} is a classically consistent system such that {\bf B} $\subseteq$ {\bf S} $\subseteq$ {\bf I} then {\bf S}$^{+g}$ is also classically consistent and {\bf B}$^{+g}$ $\subseteq$ {\bf S}$^{+g}$.  However, if  $\mathcal{Y}$ is the collection of {\em all} classically consistent subsystems of {\bf I} containing {\bf B}, then
$\bigcup \{\mathrm{\bf S} : \mathrm{\bf S} \in \mathcal{Y}\}$ is classically inconsistent, as the following argument (inspired by Vafeiadou's proof of Theorem 7.6.1) shows. Let Con({\bf B}) be a sentence of $\mathcal{L}$({\bf I}) expressing the statement ``{\bf B} $\not \vdash \; 0=1$.''

\subsubsection{\bf Theorem} Consider the intermediate systems {\bf S$_1$} $\equiv$ {\bf B} + $\mathrm{(WLPO \rightarrow Con({\bf B}))}$ and
 {\bf S$_2$} $\equiv$ {\bf B} + $\mathrm{(WLPO \rightarrow \neg Con({\bf B}))}$.
 \begin{enumerate}
 \item[(i)]{{\bf S$_1$} and {\bf S$_2$} are classically consistent subsystems of {\bf I}.}
 \item[(ii)]{{\bf S$_1$} + {\bf S$_2$} is not classically consistent.}
 \item[(iii)]{({\bf S$_1$})$^g$ conflicts with ({\bf S$_2$})$^g$, so $\bigcup \{\mathrm{\bf S}^g : \mathrm{\bf S} \in \mathcal{Y}\}$ is inconsistent.}
 \end{enumerate}

{\em Proofs.} (i) is a consequence of G\"odel's second incompleteness theorem with the fact that {\bf I} $\vdash$ $\mathrm{\neg WLPO}$.  (ii) holds because ({\bf S$_1$})$^\circ$ $\vdash$ Con({\bf B}) and ({\bf S$_2$})$^\circ$ $\vdash$ $\mathrm{\neg}$ Con({\bf B}).  (iii) follows immediately because ({\bf S$_1$})$^g$ $\vdash$ (Con({\bf B}))$^g$ and ({\bf S$_2$})$^g$ $\vdash$ $\mathrm{\neg}$ (Con({\bf B}))$^g$. \qed

\section{Alternative varieties of constructive analysis}

\subsection{Axiomatizing the recursive model}

Troelstra and van Dalen \cite{tvd1988} propose that constructive recursive mathematics RUSS, up to and including the Kreisel-Lacombe-Shoenfield-Tsejtlin  Theorem, should be axiomatized in the language of arithmetic by {\bf CRM} = {\bf HA} + ECT$_0$ + MP$_0$, where ECT$_0$ is Troelstra's ``extended Church's Thesis'' (cf. \cite{tr1973}) and MP$_0$ is an arithmetical form of Markov's Principle.  By number-realizability, {\bf CRM} is consistent relative to its classically consistent subtheory {\bf HA} + MP$_0$, but (unlike {\bf CRM}) all of Russian recursive mathematics is consistent with classical logic. A classically sound formalization of RUSS appears to require sequence variables.

The $\omega$-model of $\mathcal{L}$({\bf I}) in which the infinite sequences are the recursive sequences satisfies the classically sound theory {\bf MRA} $\equiv$ {\bf IRA} + CT$_1$ + MP$_1$ where the axiom
\[\mathrm{CT_1. \; \; \; \forall \alpha \exists e [\forall x \exists y T(e,x,y) ~\&~ \forall x \forall y (T(e,x,y) \rightarrow U(y) = \alpha(x))]}\]
(abbreviated $\mathrm{\forall \alpha GR(\alpha)}$), with no parameters allowed, plays the restrictive role of Church's Thesis.  CT$_1$ fails in {\bf B$^\circ$} by Lemma 9.8 of \cite{klve1965}, and is refutable in {\bf I} using Brouwer's Principle for Numbers ($^\ast$27.2 in \cite{klve1965}).   Its negative interpretation, however, is provable in {\bf MRA} and is consistent with {\bf I} (but not with {\bf B$^\circ$}).\footnote{In contrast, the negative interpretation of the continuous choice axiom CC$_{11}$ of {\bf I} is inconsistent with {\bf I} and with {\bf B$^\circ$}.}
\subsubsection{\bf Theorem.}
\begin{enumerate}
\item[(i)]{{\bf MRA$^{+g}$} = {\bf MRA}.}
\item[(ii)]{{\bf MRA} can be negatively interpreted in its subsystem {\bf IRA} + $\Sigma^0_1$-DNS$_0$ + $\mathrm{\forall \alpha \neg \neg GR(\alpha)}$.}
\item[(iii)]{{\bf I} + $\Sigma^0_1$-DNS$_0$ + $\mathrm{\forall \alpha \neg \neg GR(\alpha)}$ + VS is consistent and proves $\mathrm{\neg MP_1}$.}
\end{enumerate}
{\em Proofs.}  (i) holds by {\bf IA$_1$} + $\Sigma^0_1$-DNS$_0$ $\vdash$ $\mathrm{(\forall \alpha \neg \neg GR(\alpha) \leftrightarrow  (CT_1)}^g)$, Theorem 5.1(ii) and the easy facts that {\bf IRA} + MP$_1$ $ \vdash $ $\Sigma^0_1$-DNS$_0$ and
{\bf IA$_1$} + CT$_1$ $\vdash$ $\mathrm{\forall \alpha \neg \neg GR(\alpha)}$.
(ii) follows by the proof of Corollary 5.3(ii), and (iii) holds by classical $^G$realizability (\cite{jrm1971}) and \cite{ve1970}. \qed

\vskip 0.1cm

Kleene's formalization \cite{kl1969} of the theory of recursive functionals can be carried out in {\bf IRA} so constructive recursive mathematics should be formalizable in {\bf MRA}. The arithmetical recursive choice principle CT$_0$ which holds in {\bf CRM} suggests adding to {\bf MRA} either $\mathrm{AC_{00}^{Ar}}$ or a comprehension principle
\[\mathrm{CF_d^{Ar}. \; \; \forall x (A(x) \vee \neg A(x)) \rightarrow \exists \alpha \forall x (\alpha(x) = 0 \leftrightarrow A(x))}\]
restricted to formulas $\mathrm{A(x)}$ without free sequence variables.  {\bf MRA} + CF$\mathrm{_d^{Ar}}$ should be consistent by recursive number-realizability, and its minimal classical extension relative to the recursive model follows the pattern of Vafeiadou's Theorem 7.6.1.

\subsection{Bishop's constructive mathematical analysis} Anything that can be formalized in  Troelstra's {\bf EL} + AC$_{01}$, which has been used by Bishop constructivists, can also be formalized in the common subsystem {\bf IA$_1$} + AC$_{01}$ of {\bf I} and {\bf B$^\circ$} by \cite{gvf2012}.  By Theorem 5.1(vi),(xi) with the fact that {\bf IA$_1$} + AC$_{00}$ + $\mathrm{(\neg \neg A \rightarrow A)}$ $\vdash$ BI$_1$, Bishop's constructive analysis BISH has the same classical content as Kleene's {\bf B}.

\subsection{Afterword}
Reverse constructive mathematics establishes precise connections among mathematical theorems, function existence axioms, and logical principles over weak constructive theories based on intuitionistic logic.  The (weak and not so weak) base theories used here are classically sound subsystems of Kleene's formal system {\bf I} in \cite{klve1965}. Kreisel's {\em two} uncomplimentary reviews notwithstanding, Kleene and Vesley's book contained the first coherent treatment of Brouwer's intuitionistic analysis in ordinary mathematical language with intuitionistic logic, together with a proof of its consistency relative to its classically sound subtheory {\bf B}.

The classical contents (as expressed by the G\"odel-Gentzen negative interpretation) of classical analysis with countable choice, Bishop's constructive analysis, and Markov's recursive analysis are individually consistent with Kleene's and Vesley's versions \cite{klve1965}, \cite{ve1970} of intuitionistic analysis.

The perceived conflicts among CLASS, INT, BISH and RUSS partly reflect the ways language is used in these four varieties of mathematical practice.  Gentzen's negative interpretation enables a parallel treatment of classical and constructive mathematics by making  linguistic differences explicit, restricting the logic to be intuitionistic, and expressing classical reasoning in the negative language. The constructive cost of reconciliation can be measured precisely by computing the minimum classical extensions of classically sound theories.

Of course, those conflicts are never just a matter of linguistic interpretation or of intuitionistic versus classical logic.  They also reflect fundamentally different ideas about what constitutes an infinite sequence of natural numbers.

\section{Acknowledgements} Warm thanks to Garyfallia Vafeiadou for clarifying the connections among weak subsystems of constructive and intuitionistic analysis and for her collaboration in \cite{jrmgvf2021}, and to the anonymous reviewer of that paper who made me rethink this one.  Thanks again to Robert Solovay for the clever proof that suggested these questions. Many thanks to Wim Veldman for keeping the subject of intuitionistic analysis very much alive and for not allowing Kleene's pioneering work to be forgotten, and to Anne Troelstra whose mathematical legacy is secure. He is greatly missed.\footnote{This paper builds on the extended abstract for my contributed talk, dedicated to the memory of A. S. Troelstra, for the Twelfth Panhellenic Logic Symposium in Crete in June, 2019. I also thank the organizers of Computability in Europe 2021 for inviting me to talk about this subject.}

\bibliographystyle{plain}
\bibliography{calibratingthenegativeinterpretationv4forarxiv}

\end{document}